\documentclass[11pt]{article}
\usepackage{latexsym}
\usepackage{amsfonts}
\usepackage{amsmath}
\usepackage{epsfig}

\newtheorem{thm}{Theorem} \newtheorem{lemma}[thm]{Lemma}
 \newtheorem{cor}[thm]{Corollary}

\newcommand{\prob}{\mbox{\bf P}}
\newcommand{\E}{{\bf E\,}}

\newcommand{\R}{\mathbb R}

\newcommand{\eps}{c}
\newcommand{\be}{\begin{equation}}
\newcommand{\ee}{\end{equation}}
\newcommand{\C}{{\cal{C}}}
\newcommand{\F}{{\cal{F}}}

\def\qed{\relax\ifmmode\hskip2em \Box\else\unskip\nobreak\hfill
$\Box$\fi}

\newcounter{mycount}

\title{The critical random graph, with martingales}
\author{{\sc Asaf Nachmias and Yuval Peres\thanks{U.C. Berkeley
and Microsoft Research. Research of both authors
supported in part by NSF grants \#DMS-0244479 and \#DMS-0104073}}}
\begin{document}

\maketitle

\begin{abstract}
We give a short proof that the largest component $\C_1$ of the
random graph $G(n, 1/n)$ is of size approximately $n^{2/3}$. The
proof gives explicit bounds for the probability that the ratio is
very large or very small. In particular, the probability that
$n^{-2/3}|\C_1|$ exceeds $A$ is at most $e^{-cA^3}$ for some
$c>0$.
\end{abstract}

\section{Introduction}

The random graph $G(n,p)$ is obtained from the complete graph on
$n$ vertices, by independently retaining each edge with
probability $p$ and deleting it with probability $1-p$.
 Erd\H{o}s and R\'enyi~\cite{ER} introduced this model in 1960, and discovered that
as $c$ grows, $G(n,c/n)$ exhibits a {\em double jump\/}:
the cardinality of the largest component $\C_1$ is of order $\log n$
 for $c<1$,  of order $n^{2/3}$ for $c=1$ and linear in $n$ for $c>1$.
 In fact, for the critical case  $c=1$ the argument in
\cite{ER} only established the lower bound on $\prob ( | \C_1 | >
An^{2/3})$ for some constant $A>0$; the upper bound was proved
much later in \cite{B1} and \cite{L}.

Short proofs of the results stated above for the noncritical cases
$c<1$ and $c>1$ can be found  in the books \cite{AS}, \cite{B2},
and \cite{JLR}. However, we could not find in the literature
 a short and self-contained analysis of the case $c=1$.
We prove the following two theorems:

\begin{thm} [see \cite{P} and \cite{SS} for similar estimates] \label{upper}
Let ${\cal C}_1$ denote the largest component of  $G(n,1/n)$, and
let $C(v)$ be the component that contains a vertex $v$. For any
$n>1000$ and $A>8$ we have

$$ \prob ( |C(v)| > An^{2/3} ) \leq
4n^{-1/3}e^{-\frac{A^2(A-4)}{32}} \, ,$$ and
$$ \prob ( |{\cal C}_1| > An^{2/3} ) \leq \frac{4}{A}
e^{-\frac{A^2(A-4)}{32}} \, .$$
\end{thm}

\begin{thm} \label{lower}
For any $0 < \delta < 1/10$ and $n>200/\delta^{3/5}$, the random
graph $G(n,1/n)$ satisfies
$$ \prob \Big ( |\C_1| <\lfloor \delta n^{2/3} \rfloor \Big ) \leq 15\delta^{3/5}  \, .$$
\end{thm}

While the estimates in these two theorems are not optimal, they
are explicit, so the theorems say something about $G(n,1/n)$ for
$n=10^9$ and not just as $n \to \infty$. The theorems can be extended
to the ``critical window'' $p=1/n+\lambda n^{-4/3}$, see Section
\ref{secwindow}. As noted above, Erd\H{o}s and R\'enyi~\cite{ER}
proved a version of Theorem \ref{lower}; their argument was based
on counting tree components of $G(n,1/n)$.
However, to prove Theorem \ref{upper} by a similar counting
argument requires consideration of subgraphs that are not trees.
Indeed, with such considerations, Pittel \cite{P} proves tail
bounds on $n^{-2/3}|C_1|$ that are asymptotically more precise
than Theorems \ref{upper} and \ref{lower}. For a probabilistic
approach to Theorem \ref{upper} that does not use martingales, see
Scott and Sorkin~\cite{SS}.

The systematic study of the phase transition in $G(n,p)$
around the point $p\sim 1/n$ was initiated by Bollob\'as~\cite{B1} in 1984
and an upper bound of order $n^{2/3}$ for the median (or any quantile) of $|\C_1|$
was first proved by  \L uczak~\cite{L}.
 \L uczak, Pittel and Wierman~\cite{LPW} subsequently proved
the following more precise result.
\begin{thm} [\L uczak, Pittel and Wierman 1994] \label{lpw}
Let $p={1 \over n} + \lambda n^{-4/3}$ where $\lambda \in \R$ is
fixed. Then for any integer $m>0$, the sequence $$(n^{-2/3}|\C_1|,
n^{-2/3}|\C_2|, \ldots , n^{-2/3}|\C_m|)$$ converges in
distribution to a random vector with positive components.
\end{thm}
The proofs in \cite{L}, \cite{LPW} and \cite{P} are quite
involved, and use the detailed asymptotics
from \cite{Wr}, \cite{B1} and \cite{BCM} for the number 
of graphs on $k$ vertices with $k+\ell$ edges.
Aldous~\cite{A} gave a more conceptual proof of  Theorem~\ref{lpw} using diffusion approximation,
and identified the limiting distribution in terms of excursion lengths of reflected
Brownian motion with variable drift. The argument in \cite{A} is beautiful but not
elementary, and it seems hard to
extract from it explicit estimates for specific finite $n$.
A powerful approach, that works in the more general setting of percolation
on certain finite transitive graphs, was recently developed
in~\cite{BCHSS}.  This work is based on the lace expansion,
and is quite difficult.


Our proofs of Theorems \ref{upper} and \ref{lower} use an exploration process introduced in
\cite{M} and \cite{K}, and the following  classical theorem (see, e.g. \cite{D} section 4, or
\cite{W}).

\begin{thm} [Optional stopping theorem] \label{optstop}
Let $\{X_t\}_{t \ge 0}$ be a martingale for the increasing $\sigma$-fields
$\{\F_t\}$ and suppose that $\tau_1, \tau_2$ are stopping times with $0 \le \tau_1 \le \tau_2$.
If the process $\{ X_{t \wedge \tau_2} \}_{t \ge 0}$ is  bounded, then $\E X_{\tau_1} =
\E X_{\tau _2}$.
\end{thm}
\noindent{\bf Remark.} If $\{X_t\}_{t \ge 0}$ is a submartingale
(supermartingale), then under under the same boundedness
condition, we have $\E X_{\tau_1} \leq \E X_{\tau_2}$
(respectively, $\E X_{\tau_1} \geq \E X_{\tau _2}$).
\newline

The rest of the paper is organized as follows. In Section
\ref{secexplore} we present the exploration process mentioned
above. In Section \ref{seceasy} we present a very simple proof of
the fact that in $G(n,1/n)$ we have $\prob (|{\cal C}_1|>
An^{2/3}) \leq 6A^{-3/2}$.
%
The proof of Theorem \ref{upper} and
\ref{lower} are then presented in Sections \ref{secupper} and
\ref{seclower}. The technical modifications required to
handle the ``critical window'' $p=1/n+\lambda n^{-4/3}$
are presented in Section \ref{secwindow}.


\section{The exploration process} \label{secexplore}

For a vertex $v$, let $\C(v)$ denote the connected component that
contains $v$. We recall an exploration process, developed
independently by Martin-L\"of \cite{M} and Karp \cite{K}.
 In this process, vertices will be either
{\em active, explored} or {\em neutral}. At each time $t\in \{0,1, \ldots,n\}$,
the number of active vertices will be denoted $Y_t$ and
the number of explored vertices will be $t$.
Fix an ordering of the vertices, with $v$ first.
At time $t=0$, the vertex $v$ is active and all other vertices are
neutral, so $Y_0=1$.  In step $t \in \{1,\ldots,n\}$, if $Y_{t-1}>0$ let $w_t$ be the first active vertex;
if $Y_{t-1}=0$, let $w_t$ be the first neutral vertex. Denote by $\eta_t$ the number of neutral
neighbors of $w_t$ in $G(n,1/n)$, and change the status of these vertices to {\em active}. Then, set $w_t$
itself {\em explored}.

Write $N_t=n-Y_t-t-{\mathbf 1}_{\{Y_t=0\}}$.
Given $Y_1,\ldots, Y_{t-1}$,  the random variable $\eta_t$ is
distributed ${\rm Bin}(N_{t-1}, 1/n)$, and we have the recursion
\be \label{recur}
 Y_t= \left \{
\begin{array}{ll}
Y_{t-1} + \eta_t-1, & Y_{t-1} > 0 \\
\eta_t, & Y_{t-1} = 0 \, . \\
\end{array} \right .
\ee
At  time $\tau = \min \{t \ge 1: Y_t = 0\}$ the set of
explored vertices is precisely  $\C(v)$, so  $|\C(v)| = \tau$.

To prove Theorem \ref{upper}, we will couple $\{Y_t\}$ to a random walk with shifted binomial increments.
We will need the following lemma concerning the {\em overshoots\/} of such walks.
\begin{lemma} \label{l.bin} Let $p\in (0,1)$ and $\{\xi_i\}_{i \ge 1}$ be i.i.d.\ random variables with
 {\rm Bin}$(n,p)$ distribution
and let $S_t=1+\sum_{i=1}^t (\xi_i-1)$. Fix an integer $H>0$, and define 
$$ \gamma = \min \{ t \ge 1 \, :  \, S_t \geq H {\rm \ or \ } S_t = 0\} \, . $$
Let $\Xi \subset {\mathbb N}$ be a set of positive integers. Given
the event $\{S_{\gamma} \ge H , \, \gamma \in \Xi \}$, the
conditional distribution of the overshoot $S_{\gamma} - H$ is
stochastically dominated by the binomial distribution {\rm
Bin}$(n,p)$.
\end{lemma}

\noindent {\bf Proof.} First observe that if $\xi$ has  a
Bin$(n,p)$ distribution, then the conditional distribution of
$\xi-r$ given $\xi \geq r$ is stochastically dominated by
Bin$(n,p)$. To see this, write $\xi$ as a sum of $n$ indicator
random variables $\{I_j\}_{j=1}^n$ and let $J$ be the minimal
index such that $\sum _{j=1}^J I_j = r$. Given $J$, the
conditional distribution of $\xi-r$ is  Bin$(n-J,p)$ which is
certainly dominated by  Bin$(n,p)$.

For any $\ell \in \Xi $, conditioned on $ \{\gamma=\ell\} \cap
\{S_{\ell -1} = H -r\} \cap \{S_{\gamma} \ge H\}, $ the overshoot
$S_{\gamma} - H$ equals $\xi_\ell-r$ where $\xi_\ell$ has a
Bin$(n,p)$ distribution conditioned on $\xi_\ell \geq r$. The
assertion of the lemma follows by averaging. \qed

\begin{cor} \label{c.bin}
Let $X$ be distributed {\rm Bin}$(n,p)$ and let $f$ be an
increasing real function. With the notation of the previous lemma,
we have
$$ \E [ f(S_{\gamma}-H) \mid S_{\gamma} \ge H, \, \gamma \in \Xi ] \leq \E f(X) \, .$$
\end{cor}

\section{An Easy Upper Bound} \label{seceasy}

Fix a vertex $v$. To analyze the component of $v$ in $G(n,1/n)$,
we use the notation established in the previous section. We can
couple the sequence $\{\eta_t\}_{t \ge 1}$ constructed there, to a
sequence $\{\xi_t\}_{t \ge 1} $  of i.i.d.~Bin$(n,1/n)$ random
variables, such that $\xi_t \ge \eta_t$ for all $t \le n$.  The
random walk $\{S_t\}$ defined in Lemma \ref{l.bin} satisfies $S_t
= S_{t-1}+ \xi_t -1$ for all $t \ge 1$ and  $S_0=1$. Fix an
integer $H>0$ and define $\gamma$ as in Lemma \ref{l.bin}. Couple
$S_t$ and $Y_t$ such that $S_t \ge Y_t$ for all $t\leq \gamma$.
Since $\{S_t\}$ is a martingale, optional stopping gives $ 1 =
\E[S_\gamma]
 \geq H \prob (S_{\gamma}\ge H) \, ,$ 
whence
\begin{equation} \label{bound2}
\prob (S_{\gamma} \ge H) \leq {1 \over H} \, .
\end{equation}
Write $S_{\gamma} ^2 = H^2 + 2H(S_{\gamma} - H) + (S_{\gamma} -
H)^2 \, $ and apply Corollary \ref{c.bin} with $f(x)=2Hx+x^2$ to
get for $H\geq 2$ that
\begin{equation} \label{forpart2} \E \Big [ S_{\gamma}^2 \mid
S_{\gamma} \ge H \Big ] \leq H^2 + 2H + 2 \leq H^2 + 3H \, .
\end{equation}
Now $S_t^2 -(1-\frac{1}{n}) t$ is also a martingale. By optional
stopping, (\ref{bound2}) and  (\ref{forpart2}),
$$1 + (1-\frac{1}{n})\E \gamma = \E(S_\gamma^2)=\prob (S_{\gamma}\ge H) \E \Big [ S_{\gamma}^2 \mid S_{\gamma} \ge H \Big
]\leq H + 3 \, ,$$ hence we have for $2 \leq H\leq n-3$ that
\begin{equation} \label{bound3}
\E \gamma \leq H+3 \, .
\end{equation}

We conclude that for $2 \leq H\leq n-3$
$$ \prob (\gamma \geq H^2) \leq {H+3 \over H^2} \leq {2 \over H} \, .$$
Define $\gamma^* = \gamma \wedge H^2$, and so by the previous
inequality and (\ref{bound2}) we have

\be \label{bd1} \prob ( S_{\gamma^*} > 0 ) \leq \prob (S_\gamma
\geq H) + \prob ( \gamma \geq H^2 ) \leq {3 \over H} \, . \ee

Let $T=H^2$ and note that if $|C(v)| > H^2$ we must have
$S_{\gamma^*} > 0$ so by (\ref{bd1}) we deduce $\prob (|C(v)|
> T) \leq {3 \over \sqrt{T}}$ for all $ 9 \leq T \leq  (n-3)^2$. Denote by $N_{T}$
the number of vertices contained in components larger than $T$.
Then
\begin{eqnarray*} \prob \Big ( |{\cal C}_1| > T \Big ) &\leq& \prob
\Big ( |N_{T}| > T \Big ) \leq \frac{\E N_{T}}{T}
 \leq \frac{n \prob(|C(v)| > T)}{T} \, .
 \end{eqnarray*}
Putting $T= \Big ( \lfloor \sqrt{An^{2/3}} \rfloor \Big )^2$ for
any $A>1$ yields
$$\prob \Big ( |{\cal C}_1| > An^{2/3} \Big ) \leq  \prob \Big ( |{\cal C}_1| > T \Big )
\leq \frac{3n}{ \Big ( \lfloor \sqrt{An^{2/3}} \rfloor \Big )^3 }
\leq {6 \over A^{3/2}} \, ,$$ as $\Big ( \lfloor \sqrt{An^{2/3}}
\rfloor \Big )^3 \geq \Big ( \sqrt{An^{2/3}} - 1 \Big )^3 \geq
nA^{3/2}(1-3A^{-1/2}n^{-1/3}) \geq { A^{3/2}n \over 2}$. \qed

\section{Proof of Theorem \ref{upper}} \label{secupper}
We proceed from (\ref{bd1}). Define the process $\{Z_t\}$ by
\begin{equation} \label{defZ} Z_t = \sum_{j=1}^t (\eta_{\gamma^*+j} -1) \, .
\end{equation}
The law of $\eta_{\gamma^*+j}$ is stochastically dominated by a
Bin$(n-j,{1\over n})$ distribution, for $j \leq n$. Hence,
$$ \E \Big [ e^{\eps (\eta_{\gamma^*+j}-1)} \mid \gamma^* \Big ]  \leq  e^{-\eps} \Big [ 1+ {1\over n}
(e^\eps -1) \Big ]^{n-j} \leq e^{(\eps + \eps^2)\frac{n-j}{n} -
\eps} \leq e^{\eps^2 - {\eps j \over n}} \, ,$$ as $e^\eps -1 \leq
\eps + \eps^2$ for any $\eps \in (0,1)$ and $1+x \leq e^x$ for
$x\geq 0$. Since this bound is uniform in $S_{\gamma^*}$ and
$\gamma^*$, we have
$$ \E \Big [ e^{\eps Z_t} \mid S_{\gamma^*} \Big ] \leq
e^{t\eps^2 - {\eps t^2 \over 2n}} \, .$$ Write $\prob _S$ for the
conditional probability given $S_{\gamma^*}$. Then for any $\eps
\in (0,1)$, we have
$$ \prob _S \Big ( Z_t \geq -S_{\gamma*} \Big ) \leq \prob _S \Big ( e^{\eps
Z_t} \geq e^{-\eps S_{\gamma^*}} \Big  ) \leq e^{t\eps^2 - {\eps
t^2 \over 2n}} e^{\eps S_{\gamma^*}} \, .$$ By (\ref{recur}), if
$Y_{\gamma^*+j} > 0$ for all $0 \leq j \leq t-1$, then $Z_j =
Y_{\gamma^*+j} - Y_{\gamma^*}$ for all $1 \leq j\leq t$. It
follows that
\begin{eqnarray} \label{midstep}
\prob \Big ( \forall j \leq t \quad Y_{\gamma^* +j} > 0 \mid
S_{\gamma^*}>0 \Big ) &\leq& \E \Big [ \prob _S ( Z_t \geq
-S_{\gamma*}) \mid S_{\gamma*} > 0 \Big ] \nonumber \\ &\leq&
e^{t\eps^2 - {\eps t^2 \over 2n}} \E [ e^{\eps S_{\gamma^*}} \mid
S_{\gamma^*} > 0 ] \, .
\end{eqnarray}
By Corollary \ref{c.bin} with $\Xi = \{1,\ldots, H^2\}$, we have
that for $c\in (0,1)$,
\begin{equation} \label{cor6ap} \E [ e^{\eps S_{\gamma^*}} \mid \gamma \leq H^2 , \,   S_\gamma >
0] \leq e^{\eps H + \eps + \eps^2}\, .
\end{equation}
Since $\{ S_{\gamma^*} > 0 \} = \{ \gamma > H^2 \} \cup \{ \gamma
\leq H^2, \, S_\gamma > 0\}$ (a disjoint union), the conditional
expectation $\E [ e^{\eps S_{\gamma^*}} \mid S_{\gamma^*}> 0 ]$ is
a weighted average of the conditional expectation in
(\ref{cor6ap}) and  of $\E [e^{\eps S_{\gamma^*}} \mid \gamma >
H^2] \leq e^{\eps H}$. Therefore $E [ e^{\eps S_{\gamma^*}} \mid
S_{\gamma^*}> 0] \leq e^{\eps H + \eps + \eps^2}$, whence by
(\ref{midstep}), \be \label{bd2} \prob \Big ( \forall j \leq t
\quad Y_{\gamma^* + j} > 0 \mid S_{\gamma^*}>0 \Big ) \leq
e^{t\eps^2 - {\eps t^2 \over 2n} + \eps H + \eps + \eps^2} \, .
\ee By our coupling, for any integer $T>H^2$, if $|C(v)|>T$ then
we must have $S_{\gamma^*} > 0$ and $Y_{\gamma^* +j} > 0$ for all
$j \in [0,T-H^2]$. Thus, by (\ref{bd1}) and (\ref{bd2}), we have
\begin{eqnarray} \label{parabola} \prob (
|C(v)|
> T ) &\leq& \prob ( S_{\gamma^*}>0 ) \prob \Big (\forall j \in [0, T-H^2] \quad  Y_{\gamma^*+j}
> 0\mid S_{\gamma^*}>0 \Big ) \nonumber \\ &\leq& {3 \over H}
e^{(T-H^2)\eps^2 - {\eps (T-H^2)^2 \over 2n} + \eps H + \eps +
\eps^2} \, .\end{eqnarray} Take $H=\lfloor n^{1/3} \rfloor$ and
$T=\lfloor An^{2/3} \rfloor$ for some $A>4$; substituting $\eps$
which attains the minimum of the parabola in the exponent of the
right hand side of (\ref{parabola}) gives
\begin{eqnarray*} \prob ( |C(v)| > An^{2/3} ) &\leq&
4n^{-1/3}e^{-{\Big ( (T-H^2)^2/(2n) - H - 1 \Big)^2 \over
4(T-H^2+1) }} \\ &\leq& 4n^{-1/3} e^{- { \Big ( (A-1-n^{-2/3})^2/2
- 1 - n^{-1/3} \Big)^2 \over 4(A-1 + 2n^{-1/3} + n^{-2/3}) }} \leq
4n^{-1/3} e^{ - {\Big ({(A-2)^2\over 2}-2\Big )^2 \over 4(A-1/2)}}
\, ,
\end{eqnarray*}
since $H^2 \geq n^{2/3}(1-2n^{-1/3})$ and $n>1000$. As $[
(A-2)^2/2-2]^2 = A^2(A/2-2)^2$ and $(A/2-2)/(A-1/2)>1/4$ for $A>8$
we get
$$ \prob ( |C(v)| > An^{2/3} ) \leq 4n^{-1/3} e^{-{A^2(A-4) \over
32}} \, .$$


Denote by $N_{T}$ the number of vertices contained in components
larger than $T$. Then
\begin{eqnarray*} \prob \Big ( |{\cal C}_1| > T \Big ) &\leq& \prob
\Big ( |N_{T}| >T \Big ) \leq \frac{\E N_{T}}{T}
 \leq \frac{n \prob(|C(v)| >T)}{T} \, ,
 \end{eqnarray*}
and we conclude that for all $A>8$ and $n>1000$,
$$ \prob ( |{\cal C}_1| > An^{2/3} ) \leq \frac{4}{A}
e^{-\frac{A^2(A-4)}{32}} \, .$$ \qed \newline


\section{Proof of Theorem \ref{lower}} \label{seclower}
Let  $h, T_1$ and $T_2$ be positive integers, to be specified later.
The proof is divided in two stages.
In the first, we ensure, with high probability,
ascent of $\{Y_t\}$ to height $h$ by time $T_1$.
In the second stage we show that  ${Y_t}$ is likely to remain positive for $T_2$ steps.
\medskip

\noindent{\bf Stage 1: Ascent to height $h$.}
Define
$$\tau_h= \min \{t  \le T_1 \, : \, Y_t \geq h\} \, $$
if this set is nonempty, and $\tau_h=T_1$ otherwise. If
$Y_{t-1}>0$, then $Y_t^2 -Y_{t-1}^2 =(\eta_t-1)^2+
2(\eta_t-1)Y_{t-1}$. Recall that $\eta_t$ is distributed as
Bin$(N_{t-1},1/n)$ conditioned on $Y_{t-1}$, and hence if we also
require $Y_{t-1} \leq h$ then
$$
\E \Big [Y_t^2 - Y_{t-1}^2 \mid Y_{t-1} \Big ] \ge \frac{n-t-h}{n}(1-\frac{1}{n}) -2\frac{t+h}{n}\, h  \,.
$$
Next, we require that $h < {\sqrt{n}/4}$ and $t \le T_1= \lceil
\frac{n}{8h} \rceil$, whence \be \label{equ1} \E \Big [Y_t^2 -
Y_{t-1}^2 \, \Big| \,Y_{t-1} \Big ] \ge \frac{1}{2} \ee as long as
$0 < Y_{t-1} \leq h$. Similarly, (\ref{equ1}) holds if
$Y_{t-1}=0$. Thus $Y_{t \wedge \tau_h}^2 -(t \wedge \tau_h)/2 $ is
a submartingale. The proof of Lemma \ref{l.bin} implies that
conditional on $Y_{\tau_h} \ge h$, the overshoot $Y_{\tau_h}-h$ is
stochastically dominated by a Bin$(n,1/n)$ variable. So, apply
Corollary \ref{c.bin} as in (\ref{forpart2}) with $f(x)=2hx+x^2$
to get that $\displaystyle  \E Y_{\tau_h}^2 \leq h^2 + 3h \le 2h^2
\,$ for $h\geq 3$. By optional stopping,
$$ 2h^2 \geq \E Y_{\tau _h}^2
\geq \frac{1}{2}\E \tau_h  \geq \frac{T_1}{2} \prob \Bigl( \tau_h = T_1 \Big )
\, ,$$ so
\begin{equation} \label{step1}
\prob \Big (\tau_h = T_1 \Big ) \leq \frac{4h^2}{T_1} \le \frac{32h^3}{n}\, .
\end{equation}

\medskip

\begin{figure} \label{thm2}
\centering \epsfig{file=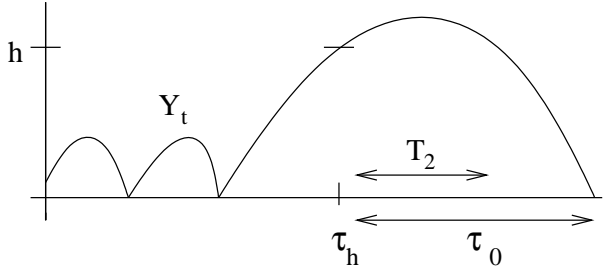} \caption{$\tau _0 \geq T_2$.}
\end{figure}

\noindent{\bf Stage 2: Keeping $Y_t$ positive for $T_2$ steps.}
\newline Define $ \tau_0 = \min \{s : Y_{\tau_h+s} =0 \}
\, $ if this set is nonempty, and $\tau_0=T_2$ otherwise. Let $M_s
= h- \min\{h,Y_{\tau_h + s}\}$. If $0<M_{s-1} <h$, then
$$
M_s^2 -M_{s-1}^2 \leq (\eta_{\tau_h+s}-1)^2
+2(1-\eta_{\tau_h+s})M_{s-1} \, ,
$$
so provided $h < \frac{\sqrt{n}}{4}$ and $s \le T_2 \le
\frac{n}{8h}$, and recalling that $\tau_h \leq T_1= \lceil
\frac{n}{8h} \rceil$ we have $ \E \Bigl[M_s^2 -M_{s-1}^2 \mid
Y_{\tau_h + s-1}, \tau_h \Bigr] \le~2 \, . $ This also holds if
$Y_{\tau_h + s -1} \ge h$, so $\{M_{S \wedge \tau_0}^2 -2(s\wedge
\tau_0)\} _{s=0}^{T_2}$ is a supermartingale. Write $\prob_h$ for
conditional probability given the event $\{Y_{\tau_h} \ge h\}$ and
$\E_h$ for conditional expectation given that event. Since $\{M_{s
\wedge \tau_0}^2 -2(s\wedge \tau_0)\}_{s=0}^{T_2}$ is a
supermartingale beginning at $0$ under $\E _h$, optional stopping
yields


\be \label{mbound} \E_h M_{\tau_0 \wedge T_2}^2 \le 2\E_h [ \tau
_0 \wedge T_2] \le 2T_2 \,, \ee whence \be \prob_h \Big ( \tau_{0}
< T_2 \Big ) \le \prob_h \Bigl(M_{\tau_0 \wedge T_2} \ge h \Bigr)
\le \frac{\E_h M_{\tau_0 \wedge T_2}^2}{h^2} \le \frac{2T_2}{h^2}
\,. \ee In conjunction with (\ref{step1}), this yields
\begin{equation} \label{step2}
\prob \Big ( \tau_0 < T_2 \Big ) \leq \prob \Big(\tau_h=T_1\Big) +
\E \prob_h \Big ( \tau_0 < T_2 \Big )
\le \frac{32h^3}{n} + \frac{2T_2}{h^2}\, .
\end{equation}
Let $T_2 = \lfloor \delta n^{2/3} \rfloor$ and choose $h$ to
approximately minimize the right-hand side of (\ref{step2}). This
gives $h = \lfloor {\delta ^{1/5} n^{1/3} \over (24)^{1/5}}
\rfloor$, which satisfies $T_2 \leq {n \over 8h}$ and makes the
right-hand side of (\ref{step2}) less than $15\delta^{3/5}$. Since
$|\C_1| <T_2$
implies $\tau_0 <T_2$, this concludes the proof. \qed \\

\section{The Critical Window} \label{secwindow}

\noindent As noted in the introduction, the proofs of Theorems
\ref{upper} and \ref{lower} can be extended to the critical
''window'' $p={1 + \lambda n^{-1/3} \over n}$ for some constant
$\lambda$. For Theorem \ref{lower} this adaptation is
straightforward, and we omit it. However, our proof of  Theorem
\ref{upper} used the fact that for $\lambda=0$ (that is, $p=1/n$)
the exploration process is stochastically dominated by a mean zero
random walk, so we include the necessary adaptation below.

\begin{thm} \label{window} Set $p={1 + \lambda n^{-1/3} \over
n}$ for some $\lambda \in {\mathbb R}$ and consider $G(n,p)$. For
$\lambda > 0$ and $A>2\lambda + 3$ we have that for large enough
$n$
$$ \prob (|C(v)| \geq An^{2/3}) \leq \Big ( \frac{4\lambda}{1-e^{-4\lambda}} +16 \Big
)n^{-1/3} e^{ - \frac{((A-1)^2/2 - (A-1)\lambda -2 )^2}{4A}} \,
,$$ and
$$ \prob ( |{\cal C}_1| \geq An^{2/3} ) \leq \Big (
\frac{4\lambda}{A(1-e^{-4\lambda})} + {16 \over A} \Big ) e^{ -
\frac{((A-1)^2/2 - (A-1)\lambda -2 )^2}{4A}} \, .$$ For $\lambda <
0$ and $A>3$ we have that for large enough $n$
$$ \prob (|C(v)| \geq An^{2/3}) \leq \Big ( \frac{-2\lambda}{e^{-\lambda}-1} + \min (5,-{1\over \lambda} ) \Big
) n^{-1/3} e^{ - \frac{((A-1)^2/2 - (A-1)\lambda -2 )^2}{4A}} \,
,$$ and
$$ \prob ( |{\cal C}_1| \geq An^{2/3} ) \leq \Big ( \frac{-2\lambda}{A(e^{-\lambda}-1)} + \min (5,-{1\over \lambda} ) \Big
) e^{ - \frac{((A-1)^2/2 - (A-1)\lambda -2 )^2}{4A}} \, .$$
\end{thm}

\noindent {\bf Proof.}  Assume $p=1/n+\lambda n^{-4/3}$ and that
$n$ is large enough; again we bound the exploration process with a
process $\{S_t\}$ defined by $S_t = S_{t-1} + \xi_t -1$ where
$\xi_t$ are distributed as Bin$(n,p)$ and $S_0 = 1$. The two cases
of $\lambda$ being positive or negative are dealt with separately;
assume first $\lambda > 0$. Since $1-e^{-a}\leq a-a^2/3$ for small
enough $a>0$, we have
$$ \E e^{-a(\xi_t -1)} = e^a[1-p(1-e^{-a})]^n \geq e^a(1-p(a-a^2/3))^n  \, .$$
By Taylor expansion of $\log(1-x)$, for small $a$ we have
\begin{eqnarray*} \log \E
e^{-a(\xi_t -1)} &\geq& a + n\Big (-p(a-a^2/3)+ O(n^{-2}) \Big ) \\
&=& a - (1+\lambda n^{-1/3})(a-a^2/3) + O(n^{-1}) \, ,
\end{eqnarray*}
and so for $a=4\lambda n^{-1/3}$ and $n$ large, we have $\E
e^{-a(\xi_t -1)} \geq 1$ hence $\{e^{-a S_t}\}$ is a
submartingale. Take $H=\lceil n^{1/3} \rceil$, and define $\gamma$
as in Lemma \ref{l.bin}. Then by optional stopping we have
$$ e^{-a} \leq 1 - \prob (S_\gamma \geq H) + \prob (S_\gamma \geq
H)e^{-aH} \, ,$$ and as $1-e^{-a} \leq a$ for $a>0$ we get

\be \label{scalebd1} \prob (S_\gamma \geq H) \leq \frac{4\lambda
n^{-1/3}}{1-e^{-4\lambda}} \, .\ee Also, observe that $S_t -
\lambda n^{-1/3} t$ is a martingale, hence by optional stopping $1
+ \lambda n^{-1/3} \E \gamma = \prob(S_\gamma \geq H) \E [S_\gamma
\mid S_\gamma \geq H]$ and so by Corollary \ref{c.bin} we get $\E
\gamma \leq \frac{8n^{1/3}}{1-e^{-4\lambda}}$. For $\lambda
> 1/4$, as $(1-e^{-4\lambda})^{-1} \leq 2$, this gives that $\E \gamma
\leq 16n^{1/3}$. It is immediate to check that $S_t^2 - {1 \over
2}t$ is a submartingale as long as $t\leq \gamma$, hence by
optional stopping ${\E \gamma \over 2} \leq \frac{4\lambda
n^{-1/3}}{1-e^{-4\lambda}} \E [S_\gamma ^2 | S_\gamma \geq H]$.
Using Corollary \ref{c.bin} as in (\ref{forpart2}) and estimating
$\frac{4x}{1-e^{-4x}} \leq 2$ for $x\in (0,1/4]$ gives the same
estimate for $\lambda \in (0,1/4]$. Thus \be \label{scalebd2} \E
\gamma \leq 16n^{1/3} \, ,\ee for all $\lambda
> 0$.  Take again $\gamma^* =
\gamma \wedge H^2$, and as in (\ref{bd1}) by (\ref{scalebd1}) and
(\ref{scalebd2}) we get \be \label{scalebd3} \prob (S_{\gamma^*} >
0) \leq \Big ( \frac{4\lambda}{1-e^{-4\lambda}} + 16 \Big
)n^{-1/3} \, . \ee

Define $Z_t$ as in (\ref{defZ}) and note that this time its
increments can be stochastically dominated by variables
distributed as Bin$(n-j,p)-1$. Similar computations to the one
made in the beginning of Section \ref{secupper} give that for
$c\in (0,1)$
$$ \E \Big [ e^{cZ_t} \mid S_{\gamma^*} \Big ] \leq e^{ ct\lambda n^{-1/3} - {ct^2 \over 2n} + c^2t(1+\lambda n^{-1/3})} \, ,$$ and so as before we
have
\begin{eqnarray*}
\prob \Big ( \forall j \leq t \quad Y_{\gamma^*+j} > 0 \mid
S_{\gamma^*}>0 \Big ) &\leq& \E \Big [ \prob _S ( Z_t \geq
-S_{\gamma*}) \mid S_{\gamma*} > 0 \Big ] \nonumber \\ &\leq& e^{
ct\lambda n^{-1/3} - {ct^2 \over 2n} + c^2t(1+\lambda n^{-1/3})}
\E [ e^{\eps S_{\gamma^*}} \mid S_{\gamma^*} > 0 ] \nonumber \\
&\leq& e^{ ct\lambda n^{-1/3} - {ct^2 \over 2n} + c^2t(1+\lambda
n^{-1/3}) + c(n^{1/3}+1) + 2(c+c^2) } \, .
\end{eqnarray*}
where the last inequality is due to Corollary \ref{c.bin}.
Write $t=\lfloor Bn^{2/3} \rfloor$ for some constant $B$ and take
$c\in(0,1)$ which attains the minimum of the parabola in the
exponent of the last display. This gives that for large enough $n$
and fixed $B>2\lambda +2$ we have
$$\prob \Big ( \forall j \leq t \quad Y_{\gamma^*+j} > 0 \mid S_{\gamma^*}>0 \Big ) \leq
e^{ - \frac{(B^2/2 - B\lambda -2 )^2}{4(B+1)}} \, . $$ Together
with (\ref{scalebd3}), as in the proof of Theorem \ref{upper}, we
conclude that for any $A>2\lambda + 3$ we have
$$ \prob (|C(v)| \geq An^{2/3}) \leq \Big ( \frac{4\lambda}{1-e^{-4\lambda}} +16 \Big
)n^{-1/3} e^{ - \frac{((A-1)^2/2 - (A-1)\lambda -2 )^2}{4A}} \,
,$$ and as before this implies that
$$ \prob ( |{\cal C}_1| \geq An^{2/3} ) \leq \Big (
\frac{4\lambda}{A(1-e^{-4\lambda})} + {16 \over A} \Big ) e^{ -
\frac{((A-1)^2/2 - (A-1)\lambda -2 )^2}{4A}} \, .$$

Assume now $p=1/n+\lambda n^{-4/3}$ for some fixed $\lambda < 0$.
For $a>0$, as $1+x \leq e^x$ we have
$$ \E e ^{a(\xi_t -1)} = e^{-a}\Big [ 1 + p(e^a-1) \Big ] ^n \leq
e^{-a+np(e^a-1)} \, .$$ By Taylor expansion of $e^x -1$ we have
$$ \log \E e ^{a(\xi_t -1)} \leq -a + (1+\lambda n^{-1/3})(a+{a^2 \over 2}+O(a^3)) \, ,$$
and so for $a = -\lambda n^{-1/3}>0$ we have that $\E e ^{a(\xi_t
-1)} \leq 1$ hence $\{e^{aS_t}\}$ is a supermartingale. With the
same $H$ and $\gamma$ as before, optional stopping gives
$$ e^a \geq 1- \prob (S_\gamma \geq H) +  \prob (S_\gamma \geq H)
e^{an^{1/3}} \, ,$$ and as $e^x -1 \leq 2x$ for $x$ small enough
we get
$$ \prob (S_\gamma \geq H) \leq {-2\lambda n^{-1/3} \over
e^{-\lambda} -1 } \, .$$ Also, as $\gamma$ is bounded above by the
hitting time of $0$, Wald's Lemma (see \cite{D}) implies that $\E
\gamma \leq -n^{1/3}/\lambda$. For $\lambda \in [-{1 \over 5} ,0]$
it is straight forward to verify that $S^2_{t \wedge \gamma} - {1
\over 2}(t \wedge \gamma)$ is a submartingale, hence as before we
deduce by optional stopping that $\E \gamma \leq 5n^{1/3}$ for
such $\lambda$'s. Thus we deduce that for all $\lambda < 0$,
$$ \E \gamma \leq \min \Big (5,-{1\over \lambda} \Big )  n^{1/3} \, .$$
The rest of the proof continues from (\ref{scalebd2}), as in the
case of $\lambda > 0$.
\qed \\

\noindent{\bf Remark.} Using similar methods, in \cite{NP}, we
analyze component sizes of bond percolation on random regular
graphs.

\section*{Acknowledgments}
We thank David Aldous, Christian Borgs and Jeff Steif for useful
suggestions and comments. Part of this research was done while the
first author was an intern at Microsoft Research. We are grateful
to the referee for a remarkably careful reading of the manuscript
and important corrections.

\bigskip \noindent
{\bf Asaf Nachmias}: \texttt{asafnach(at)math.berkeley.edu} \\
Department of Mathematics\\
UC Berkeley\\
Berkeley, CA 94720, USA.

\bigskip \noindent
{\bf Yuval Peres}: \texttt{peres(at)stat.berkeley.edu} \\
Microsoft Research\\
One Microsoft way,\\
Redmond, WA 98052-6399, USA.

\end{document}